\documentclass[a4paper,12pt]{amsart}

\usepackage[all]{xy}
\usepackage{amscd}
\usepackage{amsfonts}
\usepackage{a4wide}
\usepackage{amsmath}
\usepackage{amssymb}
\usepackage{natbib}
\usepackage{fullpage}

\numberwithin{equation}{section}

\usepackage[colorlinks=true]{hyperref}

\newcommand{\sg}{\sigma}
\newcommand{\eps}{\varepsilon}

\newcommand{\Fpbar}{\overline{\mathbb{F}}_p}
\newcommand{\ce}{\mathbb{C}}
\newcommand{\erre}{\mathbb{R}}
\newcommand{\z}{\mathbb{Z}}

\newcommand{\q}{\mathbb{Q}}

\newcommand{\h}{\mathbb{H}}

\newcommand{\XC}{X_0}

\newcommand{\G}{\Gamma_0}

\newcommand{\pruf}{\noindent \textbf{Proof}: }
\newcommand{\norm}[1]{\Vert #1 \Vert}

\newtheorem{thm}{Theorem}[section]

\newtheorem{prop-def}[thm]{Proposition-Definition}
\newtheorem{cor}[thm]{Corollary}
\theoremstyle{definition}
\newtheorem{rem}[thm]{Remark}

\title{Equidistribution of Hecke points on the supersingular module}
\author{Ricardo Menares}

\begin{document}

\begin{abstract}
For a fixed prime $p$, we consider the  (finite) set of supersingular elliptic curves over $\Fpbar$. Hecke operators act on this set. We compute the 
asymptotic frequence with which a given supersingular elliptic curve visits another under this action.
\end{abstract}

\maketitle
\section{Introduction}

Let $p$ be a prime number. We denote by $E = \{ E_1, \ldots, E_n \}$ the set of isomorphism classes of supersingular elliptic curves over $\Fpbar$. We denote by 
$S := \oplus_{i=1}^n \z E_i$ the   supersingular module in 
characteristic $p$ (i.e. $S$ is the free abelian group spanned by the elements of $E$). Hecke operators act on $S$ by
\[
T_1:= id, \quad T_m (E_i) = \sum_C E_i/C, \quad m \geq 2,
\]

\noindent where $C$ runs through the subgroup schemes of $E_i$ of rank $m$. This definition is extended by linearity to $S$ and to $S_\erre:=S\otimes \erre$. For
 an integer $m \geq 1$ we put 
\[
B_{i,j} (m) = |\{ C \subset E_i, \quad |C|=m \textrm{ and } E_i/C
\cong E_j \}|.
\]

We have that $T_m E_i = \sum_{j=1}^nB_{i,j} (m) E_j$. The  the matrix $\big(B_{i,j}(m)\big)_{i,j=1}^n$ is known as the Brandt matrix of order $m$.

For a given $D = \sum_{i=1}^n a_i E_i \in S_\erre$,
we put $\deg D = \sum_{i=1}^n a_i$.  We have that (\cite{Gross}, Proposition 2.7) $$\deg T_m E_i =  \sum_{d|m \atop p \nmid d} d=:  \sg(m)_p,$$ leading to define 
$\deg T_m:= \sg(m)_p.$

Let $M$ be the set of probability measures on $E$.  For every $i=1, \ldots , n$, we denote by $\delta_{E_i} \in M$ the Dirac measure  supported on $E_i$. Let 

\[
S^+ := \Big\{ \sum_{i=1}^n a_i E_i \in S_\erre \textrm{ such that } a_i \geq 0 \Big\} - \{0\}.
\]

For any $D = \sum_{i=1}^n a_i E_i \in S^+$, we put

\[
\Theta_D := \frac{1}{\deg D} \sum_{i=1}^n a_i \delta_{E_i}.
\]

 We have that $\Theta_D$  is a probability measure on $E$ and every element of $M$ has this form. Hence, there is a natural action of the Hecke operators on $M$, given by
$T_m \Theta_D := \Theta_{T_m D}.$

Each $E_i$ has a finite  number of automorphisms. We define

\[
w_i :=  | \textrm{Aut}(E_i)/\{\pm 1\}|, \quad W :=\sum_{i=1}^n \frac{1}{w_i} .
\]

The element $e:= \sum_{i=1}^n \frac{1}{w_i} E_i \in S \otimes \q$ is  Eisenstein (\cite{Gross}, p. 139), i.e. 
\begin{equation}\label{Eis}
T_m(e)=\deg T_m e.
\end{equation}

 We denote by  $\Theta := \Theta_e$. Equation \eqref{Eis} implies that $T_m \Theta =\Theta$ for all  $m\geq 1$. 

Let $C(E) \cong \ce^n$ be the space of complex valued functions on $E$. For $f \in C(E)$, we denote by $\norm{f}=\max_i |f(E_i)|$ and 
$$\Theta_D(f):= \int_{E} f \Theta_D = \frac{1}{\deg D} \sum_{i=1}^{n} a_if(E_i).$$ 
For a positive integer $m$, we write $m=p^km_p$ with $p\nmid m_p$. In this note, we will prove the following result:

\begin{thm} \label{principalIntro}
For all $i=1, \ldots , n$, the sequence of measures $\{ \Theta_{T_m E_i}  \}$, where $m$ runs through a set of positive integers such that $m_p$ grows to infinity, is equidistributed with 
respect to $\Theta$. More precisely, for all $\eps >0$, there exists $C_\eps >0$ such that, for every $f \in C(E)$, and for every sequence of integers $m$ such that $m_p \rightarrow \infty$, 
we have that $$|\Theta_{T_m E_i}(f) - \Theta(f)| \leq C_\eps \norm{f} n m^{-\frac{1}{2}+\eps}.$$ 
\end{thm}

 We study the asymptotic frequence of the multiplicity of  $E_j$ inside $T_m E_i$. That is, we investigate
 the behavoir of the ratio 
$B_{i,j}(m)/\deg(T_m)$ when $m$ varies. We will prove Theorem \ref{principalIntro} in the equivalent formulation:

\begin{thm} \label{principalfino}
For all $\eps>0$, there exists $C_\eps >0$ such that for every sequence of integers $m$ such that $m_p \rightarrow \infty$, we have that
\begin{equation} \label{cuantitativo}
\Big| \frac{B_{i,j}(m)}{\deg T_m} -  \frac{12}{w_j(p-1)}\Big| \leq  C_\eps m^{-\frac{1}{2}+\eps}. 
\end{equation}

 In particular, 
\begin{equation}\label{principal}
 \lim_{m_p  \rightarrow \infty} \frac{B_{i,j}(m)}{\deg T_m} = \frac{12}{w_j(p-1)}.
\end{equation}

\end{thm}

The proof of this assertion is found in section \ref{dems}.

\begin{rem}   \label{masa} The equality $\sum_{j=1}^n \frac{B_{i,j}(m)}{\deg{T_m}}  = 1,$ combined with equation \eqref{principal} implies the mass formula of Deuring and Eichler:

\[
W=\sum_{j=1}^n \frac{1}{w_j} = \frac{p-1}{12}.
\]

\end{rem}

Theorem \ref{principalIntro} can be deduced from Theorem  \ref{principalfino} as follows: remark \ref{masa} implies that 
$\Theta=\sum_{j=1}^n  \frac{12}{w_j(p-1)}\delta_{E_j}.$ Take $f \in C^0(E)$. We have that

\[
|\Theta_{T_m E_i} (f)-\Theta(f)| \leq \norm{f} \sum_{j=1}^n \Big|\frac{B_{i,j}(m)}{\deg T_m} -  \frac{12}{w_j(p-1)}\Big|.
\]

Hence, inequality \eqref{cuantitativo} implies Theorem \ref{principalIntro}.

Let $h : E \rightarrow E$ be a function. Then $h$ defines an endomorphism of $S$ and of $S_\erre$ by the rule $$h \big( \sum a_i E_i\big) := \sum a_i h(E_i).$$ 
We will also consider the  action induced on $M$ by $h^*\Theta_D:= \Theta_{h(D)}$.

\begin{cor}
 Let $q\neq p$ be a prime number.  Let $h: E \rightarrow E$ be a function such that $h \circ T_q = T_q \circ h$. Then $h^* \Theta= \Theta$. In other words, $h$ can be identified 
with a permutation $\tau \in S_n$ by $h(E_i)=E_{\tau(i)}$ and we have that $w_i = w_{\tau(i)}$ for all $i=1, \ldots, n$. 

\end{cor}

\pruf since $T_{q^k}$ is a polynomial in $T_q$, we also have that $h \circ T_{q^k} = T_{q^k} \circ h$. Let $f \in C(E)$. We have that 
\begin{eqnarray}
h^* \Theta (f) & = & \lim_{k \rightarrow \infty} h^*\Theta_{T_{q^k}E_1 } (f)  \label{pasouno} \\                    
& = & \lim_{k \rightarrow \infty}  \Theta_{h \circ T_{q^k}E_1} (f)\nonumber \\
&=&  \lim_{k \rightarrow \infty}  \Theta_{ T_{q^k} \big( h (E_1)\big)}(f) \nonumber \\
&=& \Theta (f),  \label{pasodos}
\end{eqnarray}

\noindent where we have used Theorem  \ref{principalIntro} in \eqref{pasouno} and \eqref{pasodos} $\blacksquare$ \\

The statement Theorem \ref{principalIntro}, using the Hecke invariant measure $\Theta$, has been included to emphasize the analogy with the fact that Hecke orbits
 are equidistributed on the modular curve $SL_2(\z)\backslash \h$ 
with respect to the hyperbolic measure, which is Hecke invariant (e.g. see \cite{ClozelUllmo}, Section 2).

\subsection{Weight 2 Eisenstein series for $\G(p)$}

The modular curve $\XC(p)$ has two cusps, represented by $0$ and $\infty$. We denote by $\Gamma_\infty$ (resp. $\Gamma_0$) the stabilizer of $\infty$ (resp. 0). 
The associated weight 2 Eisenstein series are given by

\begin{eqnarray}
E_\infty (z) & = & \frac{1}{2}\lim_{\varepsilon \rightarrow 0^+}
\sum_{\gamma \in \Gamma_\infty \backslash \G(p)} j_\gamma (z)^{-2}
|j_\gamma(z)|^{-2\varepsilon} \nonumber \\
E_0 (z) & = & \frac{1}{2} \lim_{\varepsilon \rightarrow 0^+}
\sum_{\gamma \in \Gamma_0 \backslash \G(p)} j_{\sigma_0^{-1}\gamma}
(z)^{-2} |j_{\sigma_0^{-1}\gamma}(z)|^{-2\varepsilon}, \nonumber
\end{eqnarray}

\noindent where $ \sigma_0 = \left(\begin{array}{cc}
0 & -1/\sqrt{p} \\
\sqrt{p} & 0
\end{array}\right)$ and $j_\eta (z) = cz+d$ for $ \eta = \left(\begin{array}{cc}
a & b \\
c & d
\end{array}\right)$ . \\

The functions $E_\infty$ and $E_0$ are weight 2 modular forms for $\G(p)$ and they are Hecke eigenforms. The Fourier expansions at  $i \infty$ are (\cite{Miyake},
Theorem 7.2.12, p. 288)

\begin{eqnarray}
E_\infty(z) & = & 1 - \frac{3}{\pi y(p+1)} + \frac{24}{p^2-1}
\sum_{n=1}^\infty b_n q^n \nonumber \\
E_0 (z) & = & - \frac{3}{\pi y (p+1)}  - \frac{24p}{p^2-1}
\sum_{n=1}^\infty a_n q^n, \nonumber
\end{eqnarray}

\noindent with the sequences $a_n$ and $b_n$ given by:

\begin{itemize}
\item if  $p \nmid n,$ then  $a_n = b_n  =  \sigma_1 (n)= \sum_{d|n}^{} d$  \\
\item if $ k \geq 1 $, then $b_{p^k}  =  p+1-p^{k+1}$ and $a_{p^k}  =  p^k $ \\
\item if $ p \nmid m$ and $k \geq 1,$  then $ b_{p^km} =
-b_{p^k}b_m \textrm{ and } a_{p^km}  =  a_{p^k} a_m$.
\end{itemize}

By taking an appropriate linear combination, we obtain a non cuspidal, holomorphic at $i\infty$ modular form

\begin{eqnarray}
f_0(z) & := & E_\infty(z) -  E_0(z) \nonumber \\
& = & 1 + \frac{24}{p^2 - 1} \sum_{n=1}^\infty (p a_n + b_n) q^n.
\nonumber
\end{eqnarray}

Since we have that

\begin{eqnarray}
E_\infty |_{\sigma_0} (z) & = & E_0(z) \nonumber \\
E_0 |_{\sigma_0} (z) & = & E_\infty (z), \nonumber
\end{eqnarray}

\noindent this shows that $f$ is holomorphic at $\G(p) 0$ as well. Since

\[
\dim_{\ce} M_2 (\G(p)) = 1 + \dim_{\ce} S_2 \big(\G(p)\big)
\]

and since $f$ is holomorphic, non zero and non cuspidal, we have the decomposition

\begin{equation} \label{sumadirecta}
M_2(\G(p)) = S_2(\G(p)) \oplus \ce f_0.
\end{equation}

\subsection{Proof of Theorem \ref{principalfino}}\label{dems}
Recall that we write $m=p^km_p$ with $p\nmid m_p$. We have that $B(p^k)$ is a permutation matrix of order dividing 2 and that $B(m)=B(p^k)B(m_p)$ (\cite{Gross}, Proposition 2.7). 
It follows that $\deg (T_m)=\deg(T_{m_p})$ and that we can define, for each $i=1, \ldots ,n$, an index $i(k) \in \{ 1, \ldots, n\}$ such that $B_{i,l} (p^k)=\delta_{i(k),l}$.   
Furthermore, $i(k)=i$ if $k$ is even.  We have that 

\begin{eqnarray}
\frac{ B_{i,j}(m)}{\deg T_m} &=&  \sum_{l=1}^n \frac{ B_{i,l}(p^k) B_{l,j}(m_p)}{\deg T_{m_p}} \nonumber \\
&=& \frac{ B_{i(k),j}(m_p) }{\deg T_{m_p}}. \nonumber
\end{eqnarray}
 
Hence, to prove Theorem \ref{principalfino} we may assume $p \nmid m$, which is what we will do in what follows. 

Our method is based on the interpretation of the multiplicities $B_{i,j}(m)$ as Fourier coefficients of a modular form.

\begin{thm} \label{teta}
For every $0 \leq i,j \leq n$, there exists a weight 2 modular form
$f_{i,j}$ for  $\G(p)$ such that its $q$-expansion at $\infty$ is
 
\[
f_{i,j}(z) := \frac{1}{2w_j} + \sum_{m=1}^\infty B_{i,j}(m) q^m,
\quad q = e^{2\pi i z}.
\]
\end{thm}

\pruf this fact is stated in \cite{Gross}, p.118. It is a particular case of \cite{Eichler}, Chapter II, Theorem 1 ($D=p, H=1, l=0$ in Eichler's notation). We remark that the theorem in \emph{loc. cit.} states modularity of a theta series constructed from an order in a quaternion algebra. The fact that this theta series is the same as our $f_{i,j}$ is a consequence of \cite{Gross}, Proposition 2.3 $\blacksquare$

Using \eqref{sumadirecta}, we can decompose

\[
f_{i,j} = g_{i,j} + c_{i,j} f_0, \quad g_{i,j} \in S_2(\G(p)), \quad c_{i,j} \in \ce.
\]

 Comparing the  $q$-expansions, we get $c_{i,j} = \frac{1}{2 w_j}.$ We have that 

$$g_{i,j}  =   f_{i,j} -c_{i,j} f_0  =  \sum_{m=1}^\infty c_m q^m,$$

\noindent where

\[
c_m = B_{i,j} (m) - \frac{12}{w_j (p^2 -1)} (p a_m + b_m).
\]

The coefficient $c_m$ depends on $(i,j)$, but we don't include this dependence in the notation in order to simplify it. Since $p \nmid m$, we have that $\deg(T_m) = \sigma_1(m)$ and

\[
c_m = B_{i,j} (m) - \frac{12}{w_j (p-1)}  \sigma_1(m).
\]

Hence,

\begin{eqnarray}
\Big|\frac{B_{i,j} (m)}{\deg T_m}  - \frac{12}{w_j (p
-1)} \Big| & = & \frac{|c_m|}{\sigma_1(m)} \nonumber \\
& \leq & \frac{|c_m|}{m}. \nonumber
\end{eqnarray}

Using Deligne's theorem (\cite{Deligne}, th\'eor\`eme 8.2, previously Ramanujan's conjecture), we have that
\[
c_m = O_\eps(m^{1/2+\eps}),
\]

\noindent concluding the proof. $\blacksquare$

\bibliographystyle{plain}

\bibliography{/home/ricardo/Escritorio/Produccion/citas}

\end{document}